\documentclass[reqno]{amsart}
\usepackage{epsfig}
\usepackage[latin1]{inputenc}
\usepackage{amssymb}
\usepackage[all]{xy}
\usepackage{color}

\newcommand{\R}{\mathrm{I\!R\!}}
\newcommand{\N}{I\hspace{-0.7ex} N}

\newcommand{\caixa}{\rule{1.3ex}{.7em}}

\newtheorem{lemma}{Lemma}[section]

\newtheorem{proposition}{Proposition}[section]

\newtheorem{theorem}{Theorem}[section]

\newtheorem{remark}{Remark}[section]

\newtheorem{ex}{Example}[section]

\begin{document}

\title[Solutions for phase field models]{On existence and uniqueness of solutions of thermal phase field models
with a general class of nonlinearities}

%\title{\Large{\bf On existence and uniqueness of solutions of thermal phase field models
%with a general class of nonlinearities}}

%\author{A. L. A. de Araujo}
%\address{Departamento de Matem\'atica, Universidade Federal de Vi\c cosa, 36571-000, Vi\c cosa (MG), Brazil}
%\curraddr{Department of Mathematics and Statistics,
%Case Western Reserve University, Cleveland, Ohio 43403}
%\email{anderson.araujo@ufv.br}
\thanks{The first author was partially supported by FAPESP, Brazil, grant 2013/22328-8.}

%\author{J.L. Boldrini}
%\address{Universidade Estadual de Campinas, Brazil.}
%\curraddr{Department of Mathematics and Statistics,
%Case Western Reserve University, Cleveland, Ohio 43403}
%\email{boldrini@ime.unicamp.br}
\thanks{The second author was partially supported by CNPq, Brazil, grant 307833/2008-9.}

%\author{B.M.R. Calsavara}
%\address{Universidade Estadual de Campinas, Brazil.}
%\curraddr{Department of Mathematics and Statistics,
%Case Western Reserve University, Cleveland, Ohio 43403}
%\email{biancamrc@yahoo.com}
\thanks{The third author was partially supported by FAPESP, Brazil, grants 2010/10087-8, 2012/15379-2.}

%    General info
\subjclass[2010]{Primary 35K61; Secundary 35Q99.}

%\date{January 1, 2001 and, in revised form, June 22, 2001.}


\keywords{Phase field model, Nonlinear parabolic systems, Leray- Schauder degree.}

\maketitle

\begin{center}
{\small A. L. A. de Araujo \\ Universidade Federal de Vi\c{c}osa, CCE, Departamento de Matem\'atica, Vi\c{c}osa, MG, Brasil \\ E-mail: \tt anderson.araujo@ufv.br}
\end{center}

\begin{center}
{\small J.L. Boldrini\\ Universidade Estadual de Campinas, IMECC, Departamento de Matem\'atica, Campinas, SP, Brasil  \\
E-mail: \tt josephbold@gmail.com}
\end{center}

\begin{center}
{\small B.M.R. Calsavara \\ Universidade Estadual de Campinas, IMECC, Departamento de Matem\'atica, Campinas, SP, Brasil  \\
E-mail: \tt biancamrc@yahoo.com}
\end{center}

%\author{\sc
%A.L.A. de Araujo
%\thanks{UFV, Brazil, anderson.araujo@ufv.br.}
%\,\,\,
%J.L. Boldrini
%\thanks{Universidade Estadual de Campinas, Brazil, boldrini@ime.unicamp.br.
%This author was partially supported by CNPq, Brazil, grant 307833/2008-9 .}
%\,\,\,
%B.M.R. Calsavara
%\thanks{Universidade Estadual de Campinas, Brazil, biancamrc@yahoo.com.
%This author was partially supported by FAPESP, Brazil, grant 2010/10087-8.}
%}

%\date{}

%\maketitle

%\thispagestyle{empty}
%\setcounter{equation}{0}
%%%%%%%%%%%%%%%%%%%%%%%%%%%%%%%%%%%%%%%%%%%%%%%%%%%%%%%%%%%%%%%%%
\begin{abstract}
We  prove the existence and uniqueness of solutions for a family of  nonlinear parabolic systems related to phase field models 
taking in account variations of temperature and the possibility of a general class of nonlinearities. 
The present results generalizes in certain aspects the already published ones in the literature.
\end{abstract}

%{\bf AMS Subject Classification:} Primary 35K61; Secundary 35Q99.

%{\bf Keywords:} Phase field model, Nonlinear parabolic systems, Leray- Schauder degree.

%%%%%%%\maketitle

%\section*{This is an unnumbered first-level section head}
%This is an example of an unnumbered first-level heading.

%% The correct journal style for \specialsection is all uppercase; a known bug
%% in amsart.cls prevents this, so input must be uppercase until it is fixed.
%\specialsection*{This is a Special Section Head}
%\specialsection*{THIS IS A SPECIAL SECTION HEAD}
%This is an example of a special section head%
%%%%%%%%%%%%%%%%%%%%%%%%%%%%%%%%%%%%%%%%%%%%%%%%%%%%%%%%%%%%%%%%%%%%%%%%
%\footnote{Here is an example of a footnote. Notice that this footnote
%text is running on so that it can stand as an example of how a footnote
%with separate paragraphs should be written.
\par
%And here is the beginning of the second paragraph.}%
%%%%%%%%%%%%%%%%%%%%%%%%%%%%%%%%%%%%%%%%%%%%%%%%%%%%%%%%%%%%%%%%%%%%%%%%

\section{Introduction}
$\indent$
Let $\Omega \subset \R^N$ be a bounded domain and $T>0$ a finite constant. We denote $Q=\Omega \times (0,T)$ and $\sum = \partial\Omega \times (0,T)$ and  consider the following nonlinear parabolic system
\begin{equation}
\label{P.1}
\left\{
\begin{array}{lll}
\displaystyle u_t + l\phi_t = \Delta u + f(x,t) &\textup{in}& Q,\\
\displaystyle  \phi_t = \Delta \phi + F(x,t,\phi) + u &\textup{in}& Q,\\
\displaystyle  \partial u/\partial \nu = \partial \phi /\partial \nu = 0 &\textup{on}& \sum,\\
\displaystyle u(x,0)=u_0(x), \ \phi(x,0)=\phi_0(x) &\textup{for all}&x \in \Omega  ,
\end{array}
\right.
\end{equation}

\noindent
that may be used to model the solidification/melting process of pure materials in the region $\Omega$ when one takes in consideration the variations of temperature.
Here the unknowns are the functions $u$ and $\phi$, related respectively to the temperature  distribution of the material and the phase function, called phase-field, used to distinguish between the liquid and solid phases of the material.
$F(x,t,\phi)$ in the phase-field equation is related to the derivative with respect to $\phi$ of the atomic interaction potential of the material being considered;
$f(x,t)$ is related to the density of heat sources or sinks, and the constant $l >0$ is related to the latent heat. Since they are not relevant for the purposes the discussion in this article, to simplify the notation, we assumed values the other material constants to be one; our results would exactly the same for other values of these constants.

We remark that the nonlinearity in~(\ref{P.1}) is given by $F(x,t,\phi)$ and that to model different realistic situations one must consider different types of such nonlinearities. Thus, it is important to understand the behavior of such systems for several classes of $F(x,t,\phi)$,  and in particular,  to know results on existence of solutions, their regularity, uniqueness, asymptotic behavior, and so on.

The first rigorous mathematical analysis for a phase-field model as in~(\ref{P.1}) was done by  Caginalp in~\cite{Caginalp} by considering the
nonlinearity coming from the classical two-wells potential: $F(\phi)=(\phi-\phi^3)/2$. Along the time other types of non-linearity were
considered, among others are the works of Hoffman and Jiang~\cite{Hoffman}, Bates and Zheng~\cite{Bates},  and Moro\c{s}anu,
Motreanu~\cite{Morosanu-Motreanu 2,Morosanu-Motreanu 1} and C\^arja~\textsl{et all}~\cite{OAC} 
(this last paper considers the case of non homogeneous boundary conditions).

\vspace{0.2cm}
Our objective in this article is to present a result on existence, uniqueness and regularity of solutions of system~(\ref{P.1}) for a class of nonlinearities $F(x,t,\phi)$ that includes the one derived from the two-wells potential in the usual tridimensional case and is related to the class considered in Moro\c{s}anu and Motreanu~\cite{Morosanu-Motreanu 2,Morosanu-Motreanu 1}. It is a bit difficult to compare in general terms ours and Moro\c{s}anu and Motreanu's classes, but, as we will show, at least for the usual and important subclass of autonomous and homogeneous nonlinearities, that is, nonlinearities that depend only on $\phi$, our subclass is strictly larger than the corresponding subclass in~\cite{Morosanu-Motreanu 1}.
Details of these comparisons and further commentaries are presented in Section~\ref{main result}.

\vspace{0.1cm}
As for the techniques we use to prove our results, similarly as in the works of Hoffman and Jiang~\cite{Hoffman} and Moro\c{s}anu and
Motreanu~\cite{Morosanu-Motreanu 2,Morosanu-Motreanu 1}, at certain point we use the Leray-Schauder fixed point theorem.

\vspace{0.1cm}
The paper is organized as follows. In Section $2$ we fix the notations, and recall certain concepts and results that will use later on; in Section $3$ we explicitly state our technical assumptions and our main result concerning existence, regularity and uniqueness of solutions; we also discuss with some detail the relations between ours and the results in~\cite{Morosanu-Motreanu 2, Morosanu-Motreanu 1}; in Section $4$, we
analyze an auxiliary problem related to the phase-field equation in~(\ref{P.1}).
We remark that the same auxiliary problem was already treated in~\cite{Morosanu-Motreanu 2, Morosanu-Motreanu 1}, but we have to reconsider it because in the present work we need some new specific arguments. Finally, Section $5$ is dedicated to the proof of our main result.

%%%%%%%%%%%%%%%%%%%%%%%%%%%%%%%%%%%%%%%%%%%%%%%%%%%%%%%%%%%%%%%%%%%%%%%%%%%%%%%%%%%%%%%%%%%%%%%%%%%%%%%%%%%%%%%%%%%%%%%%%%%%%%%%%%%%%%%%%%%%%%%%%%%%%%%%%%%%%%%%%%%%%%%%%%%%%%%%%%%%%%%%%%%%%%%%%%%%%%%%%%%%%%%%%%%%%%%%%%%%%%%%%%%%%%%%%%%%%%%%%%%%%%%%%%%%%%%%%%%%%%%%%%%%%%%%
\section{Preliminaries}
$\indent$Let $\Omega \subset \R^N$ be an open and bounded domain with a sufficiently smooth boundary, $\partial \Omega \in C^{3}$ and $Q = \Omega\times (0, T)$ the space-time cylinder with lateral surface $\sum = \partial \Omega \times (0,T)$. For $t \in (0, T]$,
we denote $Q_t = \Omega\times (0, t)$.

Given $X, Y$ Banach spaces, we denote by $X \hookrightarrow Y $ the \textsl{continuous immersion} of $X$ in $Y$ and by
$\xymatrix{\displaystyle X \ar @{^{(}->>}[r] & Y }$ the \textsl{compact immersion} of $X$ in $Y$.

\vspace{0.2cm}
Next, to ease the references, we state the following embedding result for Sobolev spaces of type $W^{r,s}_p(Q)$,
which is a particular case of Lemma $3.3$ in Ladyzhenskaya et all~\cite[pp.~80]{Ladyzhenskaya}, obtained by taking $l = 1$ and $r = s = 0$ in there.
\begin{lemma}\label{imersao continua em Lp 01}
Let $\Omega$ a domain of $\R^N$ with boundary $\partial \Omega$ (at least satisfying the cone property). Then for any function $u \in W^{2,1}_p(Q)$ we also have $u \in L^{q}(Q)$, and it is valid the following inequality
\begin{eqnarray*}\label{i.01}
\displaystyle \|u\|_{L^{q}(Q)} \leq C\|u\|_{W^{2,1}_p(Q)},
\end{eqnarray*}
provided that: $q=\infty$ if $p> \frac{N+2}{2}$; $q \geq 1$ if $p =\frac{N+2}{2}$ and $q= \frac{p(N+2)}{N+2-2p}$ if $p < \frac{N+2}{2}$. The constant $C$ depends only on $T, p, q, N$ and $\Omega$.
\end{lemma}

%%%%%%%%%%%%%%%%%%%%%%%%%%%%%%%%%%%%%%%%%%%%%%%%%%%%%%%%%%%%%%%%%%%%%%%%%%%%%%%%%%%%%%%%%%%%%%%%%%%%%%%%%%%%%%%%%%%%%%%%%%%%%%%%%%%%%%%%%%%%%%%%%%%%%%%%%%%%%%%%%%%%%%%%%%%%%%%%%%%%%%%%%%%%%%%%%%%%%%%%%%%%%%%%%%%%%%%%%%%%%%%%%%%%%%%%%%%%%%%%%%%%%%%%%%%%%%%%%%%%%%%%%%%%%%%%%%%%%%%%%%%%%%%%%%%%%%%%%%%%%%%%%%%%%%%%%%%%%%%%%%%%%%%%%%%%%%%%%%%%%%%%
\section{Main result, comparisons and commentaries}
\label{main result}

$\indent$For the rest of this article we use the following technical hypotheses:
\begin{description}
\item[$(H_0)$] $f$ (or $g) \in L^p(Q)$ with $p\geq 2$, $\phi_0$ (or $u_{0}) \in W^{2-2/p}_p(\Omega)$ such that\\
$\displaystyle \partial\phi_0/\partial\nu=0$ (or $\displaystyle  \partial u_0/\partial \nu = 0$)
on $ \sum = \partial \Omega \times (0,T)$.

\item[$(H_1)$] There is a constant $a_0 \in \R$ such that for any $(x,t)\in Q$ and $z_1,z_2 \in \R$, we have
\begin{eqnarray*}\label{}
	(F(x,t,z_1) - F(x,t,z_2))(z_1-z_2)\leq a_0(z_1-z_2)^2.
\end{eqnarray*}

\item[$(H_2)$] There is a function $G :Q\times \R^2 \rightarrow \R$ such that for any $(x,t) \in Q$ and  $z_1,z_2 \in \R$, there holds
	\[
	(F(x,t,z_1)-F(x,t,z_2))^2\leq G(x,t,z_1,z_2)(z_1-z_2)^2 , 
	\]
	
	\noindent
	where $G(x,t,z_1,z_2)\leq c_0(1+|z_1|^{2r-2}+|z_2|^{2r-2})$ with constants $c_0$ and $r \geq 1$.

\item[$(H_3)$] For $N \in \N - \{0\}$, the allowed values of the parameter  $r$ in the previous hypothesis are either
any $r \geq 1$ when $ p \geq (N+2)/2$ or $ 1 \leq r < \frac{N+2}{N +2- 2p}$ when $p < (N+2)/2$.

 \item[$(H_4)$] $F: Q\times \R \rightarrow \R$ is a Caratheodory function, i.~e., $F(.,.,z)$ is measurable on $Q$, $\forall z \in \R$, and $F(x,t,.) \in C(\R,\R)$, $\forall (x,t) \in Q$, $F(.,.,0) \in L^{\infty}(Q)$.
\end{description}

\hspace{0.3cm}
For future use in the comparison with the results of~\cite{Morosanu-Motreanu 2},~\cite{Morosanu-Motreanu 1}, and also to obtain better regularity results for the solutions, we state below the following results.

\begin{lemma}
\label{estF2}
Assumption ($H_2$) and ($H_4$) implies that $F$ fulfills the growth condition
\begin{equation}\label{estF}
	|F(x,t,z)|\leq a(1 + |z|^r),\forall (x,t,z) \in Q\times \R,
\end{equation}
where $a>0$ is a constant.
\end{lemma}

\noindent
{\bf Proof:} By setting $z_1=z$ and $z_2=0$ in $(H_2)$, we get
$$
\begin{array}{lll}
|F(x, t , z)| &\leq & |F(x, t , 0)| + G^{1/2}(x,t,z,0)|z|\\
&\leq & |F(x, t , 0)| + c_0^{1/2}(1+ |z|^{2r-2})^{1/2}|z|, \ \forall \ z \in \R.
\end{array}
$$
Since $F(.,.,0) \in L^{\infty}(Q)$, relation~(\ref{estF}) follows. \hfill\caixa

\vspace{0.1cm}
\begin{lemma}
\label{estF2.2}
Assumption ($H_1$) and ($H_4$) implies that for some constant  $d_0>0$,
\begin{equation}\label{2.2.1}
	F(x,t,z)z \leq d_0(1 + z^2), \quad \forall (x,t) \in Q, \, z \in \R.
\end{equation}
\end{lemma}

\noindent
{\bf Proof:} By setting $z_1=z$ and $z_2=0$ in $(H_1)$, we get
$F(x, t , z)z \leq  F(x, t , 0)z + a_0z^2$,  $\forall \ z \in \R$ .
Since $F(.,.,0) \in L^{\infty}(Q)$, we have
$F(x, t , z)z \leq  \|F(\cdot, \cdot , 0)\|^2_{L^{\infty}(Q)}/2 + z^2/2 + a_0z^2$, $ \forall \ z \in \R$,
and the result follows with $d_0=\max\left \{\|F(\cdot, \cdot , 0)\|^2_{L^{\infty}(Q)}/2, 1/2+a_0\right\}$. \hfill\caixa

\vspace{0.2cm}
The main purpose of the present article is to prove the following result
\begin{theorem}\label{meantheorem}
Under assumptions $(H_0)-(H_4)$. Then the problem~(\ref{P.1}) has a unique solution $(u,\phi) \in W^{2,1}_p(Q)\times W^{2,1}_p(Q)$. The solution $(u,\phi)$ satisfies
\begin{equation}\label{estprinc}
\|u\|_{W^{2,1}_p(Q)} + \|\phi\|_{W^{2,1}_p(Q)}\leq C(1+\|\phi_0\|_{W^{2-2/p}_p(\Omega)}+\|u_0\|_{W^{2-2/p}_p(\Omega)}+\|f\|_{L^p(Q)}),
\end{equation}
where the constants $C$ depends only $|\Omega|, \alpha, T, p, r, c_0, a, d_0$.
\end{theorem}

The proof of this result will be given in Section $5$.
But before proceeding with the preparation for such proof, 
in the following we will compare ours and the corresponding results in the articles by Moro\c{s}anu and Motreanu~\cite{Morosanu-Motreanu 2},~\cite{Morosanu-Motreanu 1};
just concerning the class of nonlinearities being considered, 
we also briefly compare ours with the one in the article by C\^arja~\textsl{et all}~\cite{OAC}.

\vspace{0.1cm}
In the paper of Moro\c{s}anu and Motreanu~\cite{Morosanu-Motreanu 2}, the authors used the following assumptions about the non-linearity $F$:

\vspace{0.1cm}
\begin{description}
\item[$(M_{1})$] There is a constant $a_0 \in \R$ such that for any  $(x,t)\in Q$ and  $z_1,z_2 \in \R$ there holds
$$
	(F(x,t,z_1) - F(x,t,z_2))(z_1-z_2)\leq a_0(z_1-z_2)^2. %\label{aux.7}
$$
\item[$(M_{2})$] There is a function $\overline{F}:Q \times \R^2 \rightarrow \R$ such that for any $(x,t) \in Q$ and $z_1,z_2 \in \R$, we have
$$ %\label{aux.8}
	(F(x,t,z_1)-F(x,t,z_2))^2\leq \overline{F}(x,t,z_1,z_2)(z_1-z_2)^2,
$$

\noindent
where $\overline{F}(x,t,z_1,z_2)\leq c_0(1+|z_1|^{2r-2}+|z_2|^{2r-2})$,
with some constants $c_0$ and $r \geq 1$.

\item[$(M_{3})$] There exist functions $k$ and $h$ such that $F(x,t,z)= k(x,t,z)-h(z)$, $\forall (x,t) \in Q$, $z \in \R$, where $k: Q\times \R \rightarrow \R$ is a Caratheodory function, i.~e., $k(.,.,z)$ is measurable on $Q$, $\forall z \in \R$, and $k(x,t,.) \in C(\R,\R)$, $\forall (x,t) \in Q$, $F(.,.,0) \in L^{\infty}(Q)$ and $h \in C^1(\R,\R)$. In addition, the function $k$ and $h$ verify the assumptions
\begin{description}
	 \item[$(i)$] $k(x,t,z)^2\leq a_1h(z)z + a_2(1+z^2)$, $\forall (x,t,z) \in Q\times \R$, for some constants $a_1>0$, $a_2>0$,
	 \item[$(ii)$] $-b_0 \leq h'(z)\leq b_1(1 + |z|^{r-1}) $, $\forall z \in \R$ and $1\leq r < (N+2)/N-2$ if $N>2$, for some constants $b_0, b_1>0$.
\end{description}	
\end{description}

\vspace{0.1cm}
With assumption~$(M_{2})$ the authors proved~(\ref{estF}), i.e,
\[
|k(x,t,z)-h(z)|\leq a(1 + |z|^r), \quad \forall\ (x,t,z) \in Q\times \R, \mbox{ where } a>0 \mbox{ is a constant. }
\]

\begin{lemma}\label{rel.h-F}
The condition $(M_{3})(i)$ implies the relation
\begin{eqnarray*}\label{n1.3.5}
F(x,t,z)z \leq d_0(1 + z^2), \, \, \, \,
\forall (x,t) \in Q, \, z \in \R,
\end{eqnarray*}
where $d_0$ is a positive constant. That is, the above condition $(M_{3})(i)$ implies (\ref{2.2.1}).
\end{lemma}

\noindent
{\bf Proof:} The condition~$(M_{3})(i)$ and Young's inequality show that
\begin{eqnarray*}\label{n1.5}
\begin{array}{lll}
\displaystyle
k(x,t,z)z = \left(\frac{\sqrt{2}}{\sqrt{a_1}}k(x,t,z)\right)\left(\frac{\sqrt{a_1}}{\sqrt{2}}z\right)&\leq & \displaystyle\frac{1}{a_1}k(x,t,z)^2+\frac{a_1}{4}z^2
\\
\displaystyle
&\leq & h(z)z + d_0(1 + z^2),
\end{array}
\end{eqnarray*}
for all $(x,t)\in Q$ and $z\in \R$.

Hence
$$(k(x,t,z)-h(z))z \leq d_0(1 + z^2), \, \, \, \, \,
\forall (x,t)\in Q, \, z\in \R,
$$

\noindent
and so (\ref{2.2.1}) is established. \hfill\caixa	

\vspace{0.3cm}
In the paper of Moro\c{s}anu and Motreanu~\cite{Morosanu-Motreanu 1}, the authors used on non-linearity $F$
the assumptions~(\ref{estF}), $(M_{1})$, $(M_{2})$ and

\vspace{0,2cm}

\begin{description}
  \item[$(M_4)$]  $F(x,t,z)|z|^{pr-r-1}z \leq \alpha(1 + |z|^{pr-1}) - \beta|z|^{pr}$, for constants $\alpha, \beta >0$ and $r\geq 1$ provided
$$r < \frac{N+2}{N+2 - 2p}, \; \mbox{ if } \; N+2 - 2p>0.$$
\end{description}

\begin{remark}

It follows from Lemma~\ref{rel.h-F}, $(M_{1})$ and $(M_{2})$ that the nonlinearity, $F(x,t,z)$ considered in of Moro\c{s}anu and Motreanu~\cite{Morosanu-Motreanu 2} satisfies our hypotheses $(H_{1})-(H_{4})$.

\vspace{0.1cm}
In the same way, one can prove that the class of nonlinearity considered in Moro\c{s}anu \cite{Morosanu-Motreanu 1} also satisfies our hypotheses $(H_{1})-(H_{4})$.

\vspace{0.1cm}
As for the paper of C\^arja~\textsl{et all}~\cite{OAC}, the authors consider an autonomous non-linearity $F$ of form
\[
F(\phi) = f(\phi) - a_s|\phi|^{s-1}\phi, \,\, \forall \, \phi \in \mathbb{R},
\]

\noindent
with $a_s>0$ and $s\geq 3$ satisfying
\[
s < \frac{N+2}{N+2 - 2p}, \; \mbox{ if } \; N+2 - 2p>0,
\]

\noindent
while $f(\phi) \in C^1(\mathbb{R})$ fulfills, for constants $b_1, b_2>0$, the following properties:
\[
|f'(\phi)|\leq b_1(1+|\phi|^{s-2}), \,\, \forall \, \phi \in \mathbb{R}
\]

\noindent
and
 \[
 (f(\phi_1) - f(\phi_2))(\phi_1 - \phi_2)\leq b_2(\phi_1 - \phi_2)^2, \,\, \forall \, \phi_1, \phi_2 \in \mathbb{R}.
 \]

\noindent
In \cite[Lemma 1.1]{OAC} it is proved that $F$ satisfies $(M_{1})$, $(M_{2})$ and $(M_{4})$ with $r=s$.
Therefore, by just comparing the allowed class of nonlinearities, 
 since \cite{OAC} is also concerned with non homogeneous boundary conditions, 
 the nonlinearities in \cite{OAC} also satisfies the hypotheses of the present work.

%It is enough to check hypothesis $(H_4)$.
	
	%Indeed, by $(M_{1})$ with $z_1=z$ and $z_2=0$, we have
%\[
%	(F(x,t,z) - F(x,t,0))z\leq a_0z^2,  \ \forall \ (x,t)\in Q, \ z \in \R.
%\]	
%Therefore,
%\[
%	F(x,t,z)z  \leq a_0z^2 + F(x,t,0)z \leq a_0z^2 + \frac{1}{2}|F(x,t,0)|^2 + \frac{1}{2}|z|^2,  \ \forall \ (x,t)\in Q, \ z \in \R.
%\]	
%Since $F(\cdot,\cdot,0) \in L^{\infty}(Q)$, we have the assumption $(H_4)$ verified.
\end{remark}

The previous remarks raise the possibility that our results generalize the ones of \cite{OAC}, \cite{Morosanu-Motreanu 2} and \cite{Morosanu-Motreanu 1}. This in fact is so in the sense to be described in the following.		

The results in the present work improve the results of  Moro\c{s}anu and Motreanu~\cite{Morosanu-Motreanu 2} in the sense that we have a larger range for  the allowed values of the parameter $r$ appearing in hypothesis $(H_2)$.
In~\cite{Morosanu-Motreanu 2}, it is required that $1 \leq r < (N + 2)/(N -2)$, when $N \geq 3$,
whereas here, according  to $(H_3)$ we improves to
either any $r \geq 1$ when $ p \geq (N+2)/2$ or $ 1 \leq r < (N+2)/(N +2- 2p)$ when $p < (N+2)/2$.
Observe  that this last condition  is exactly the same as in  Moro\c{s}anu Motreanu~\cite{Morosanu-Motreanu 1}, but in that paper a further restriction on the nonlinearity is required (see $(M_4)$).
In fact, in the following we will give an example of a function satisfying our hypotheses of nonlinearity but not all the assumptions of $(M_4)$ in \cite{Morosanu-Motreanu 1};  
this means that, at least in the autonomous case, our class of nonlinearities is also strictly larger than the one in \cite{Morosanu-Motreanu 1}.

\begin{ex}\label{example2}
Fix $p >2$ and consider $r_1, r_2$ satisfying
\begin{equation}
\label{in.1.2}
\begin{array}{l}
\displaystyle
1\leq r_2 < r_1 < \frac{N+2}{N+2 - 2p},   \mbox{ with } N+2 - 2p>0,
\\
\displaystyle
r_1 < r.
\end{array}
\end{equation}

Then  $F:\R \rightarrow \R$ defined by
\begin{eqnarray}\label{in.2.1}
F(z) = |z|^{r_2-1}z - |z|^{r_1-1}z, \ \ \forall \ z\in \R,
\end{eqnarray}

\noindent
satisfies the conditions of the present article but does not satisfies $(M_4)$.
\end{ex}

To check the previous claim, we start by showing that the function $F$ defined by~(\ref{in.2.1}) satisfies our hypothesis ($H_1$), ($H_2$) and ($H_4$) when
\[
r\geq r_1.
\]

From~(\ref{in.1.2}), (\ref{in.2.1}) and Mean Value Theorem we have: $\exists t_0\in (0,1)$ such that
$$
\begin{array}{l}
(F(z_1) - F(z_2))(z_1 - z_2) \\
= (r_2|(1-t_0)z_1 + t_0z_2|^{r_2-1} - r_1|(1-t_0)z_1 + t_0z_2|^{r_1-1})(z_1 - z_2)^2,	
\end{array}
$$
$\forall\ z_1$, $\ z_2 \ \in \R$.

But, the function $g(t)= r_2|t|^{r_2-1} - r_1|t|^{r_1-1}$ is such that
$$g(t)\leq 0, \ \textup{if} \ |t| \geq \left(\frac{r_2}{r_1}\right)^{\frac{1}{r_1-r_2}}$$
and
$$g(t)\leq r_2\left(\frac{r_2}{r_1}\right)^{\frac{r_2-1}{r_1-r_2}} + r_1\left(\frac{r_2}{r_1}\right)^{\frac{r_1-1}{r_1-r_2}}, \ \textup{if} \ |t| < \left(\frac{r_2}{r_1}\right)^{\frac{1}{r_1-r_2}}.$$
So there exists a constant $a_0>0$ such that
$$g(t) \leq a_0, \forall \ t \in \R.$$
Therefore
$$
\begin{array}{l}
(F(z_1) - F(z_2))(z_1 - z_2) \leq a_0(z_1 - z_2)^2, \forall z_1, \ z_2 \ \in \R,	
\end{array}
$$
which ensures that our hypothesis ($H_1$) is satisfied.

From~(\ref{in.2.1}) and Mean Value Theorem we have
$$
\begin{array}{lll}
|F(z_1) - F(z_2)| &\leq & \left| |z_1|^{r_2-1}z_1 - |z_2|^{r_2-1}z_2 \right|
+\left| |z_1|^{r_1-1}z_1 - |z_2|^{r_1-1}z_2 \right|\\
&\leq & G(z_1,z_2)|z_1 - z_2|, \ \ \forall z_1, \ z_2 \ \in \R,	
\end{array}
$$
where
$$
\begin{array}{l}
\displaystyle
G(z_1,z_2)  = r_2\sup_{0\leq t \leq 1}|(1-t)z_1 + tz_2|^{r_2-1} + r_1\displaystyle\sup_{0\leq t \leq 1}|(1-t)z_1 + tz_2|^{r_1-1}
\\
\displaystyle
\phantom{G(z_1,z_2)} \leq  r_2(|z_1| + |z_2|)^{r_2-1} + r_1(|z_1| + |z_2|)^{r_1-1}
\\
\displaystyle
\phantom{G(z_1,z_2)} \leq  a( 1 + |z_1|^{r_1-1} + |z_2|^{r_1-1})
\\
\displaystyle
\phantom{G(z_1,z_2)} \leq  a( 1 + |z_1|^{r-1} + |z_2|^{r-1}), \qquad \qquad \qquad \forall z_1, \ z_2 \ \in \R,	
\end{array}
$$
for some constant $a>0$; which ensures that our hypothesis ($H_2$) is satisfied with $F(x,t,z)= F(z)$.

Moreover, from~(\ref{in.1.2}) and~(\ref{in.2.1}) we have constants $c_1, c_2 >0$ such that
$$|F(z)| \leq |z|^{r_2} + |z|^{r_1} \leq c_1(1 + |z|^{r_1}), \ \ \forall z \ \in \R,$$
and
$$F(z)z= |z|^{r_2+1} - |z|^{r_1+1}= |z|(|z|^{r_2} - |z|^{r_1})\leq c_2(1 + z^2), \ \ \forall z \ \in \R,$$
because, $|z|^{r_2} - |z|^{r_1} \leq 1$ if $|z|\leq 1$ and by~(\ref{in.1.2}), $|z|^{r_2} - |z|^{r_1} \leq 0$ if $|z|>1$, because $1\leq r_2<r_1$. The last two inequalities ensures that our hypothesis ($H_4$) is satisfied with $F(x,t,z)= F(z)$.

\vspace{0.1cm}
Next, we will show that $F$ defined in~(\ref{in.2.1}) does not satisfy the assumption~$(M_4)$ in Moro\c{s}anu and Motreanu~\cite{Morosanu-Motreanu 1}
 (same as hypothesis~$(H_0)$ in C\^arja~\textsl{et all}~\cite{OAC}), for any $r_1< r < (N+2)/(N+2 - 2p)$.

If $1\leq r_2 < r_1 < r$, we will show that $F$ defined in~(\ref{in.2.1}) does not satisfy the inequality
$(M_4)$.

Indeed, suppose there is $\alpha, \beta$ such that $F$ satisfies $(M_4)$ for all $z \in \R$. On other hand, we have that
$$
\begin{array}{l}
\displaystyle
F(z)|z|^{pr-r-1}z =\left( |z|^{r_2-1}z -|z|^{r_1-1}z \right)|z|^{pr-r-1}z
\\
\displaystyle
\phantom{F(z)|z|^{pr-r-1}z}
=|z|^{pr+r_{2}-r} -|z|^{pr+r_{2}-r}
\end{array}
$$
Therefore,
$$|z|^{pr+r_2-r} - |z|^{pr+r_1 - r}\leq  \alpha(1 + |z|^{pr-1}) - \beta|z|^{pr},$$
which can be written as
$$|z|^{pr+r_2-r} + \beta|z|^{pr}\leq  \alpha(1 + |z|^{pr-1}) + |z|^{pr+r_1-r}.$$
Dividing this last inequality by $|z|^{pr}$, with $z \neq 0$, we obtain that
$$|z|^{r_2-r} + \beta \leq \alpha\left[\frac{1}{|z|^{pr}} + \frac{1}{|z|}\right] + |z|^{r_1-r},$$
for all $z \in \R - \{0\}$. This is a contradiction because,
$$\lim_{|z| \rightarrow +\infty} |z|^{r_2-r} +\beta = \beta>0$$
and	
$$\lim_{|z| \rightarrow +\infty} \alpha\left[\frac{1}{|z|^{pr}} + \frac{1}{|z|}\right] + |z|^{r_1-r} = 0,$$
because $1\leq r_2 < r_1 < r$.

\begin{remark}
If we restrict ourselves to the subclass of autonomous nonlinearities, i.e., nonlinearities of form $F(x,t,z) = F(z)$, our last example shows that our hypotheses $(H_{1})-(H_{4})$ are less restrictive than the hypotheses~$(M_{1})-(M_{4})$ in Moro\c{s}anu and Motreanu~\cite{Morosanu-Motreanu 1}. By \cite[Lemma 1.1]{OAC} our hypotheses $(H_{1})-(H_{4})$ are less restrictive than the hypotheses \cite{OAC}.
\end{remark}

\begin{remark}
By taking $N=3$, $k(x,t,z)=a(x,t)z + b(x,t)z^2$, with $a,b \in L^{\infty}(Q)$ and $h(z)=z^3$, problem~(\ref{P.1})  generalizes the parabolic problem studied by Hoffman and Jiang in~\cite{Hoffman}.

\end{remark}

%%%%%%%%%%%%%%%%%%%%%%%%%%%%%%%%%%%%%%%%%%%%%%%%%%%%%%%%%%%%%%%%%%%%%%%%%%%%%%%%%%%%%%%%%%%%%%%%%%%%%%%%%%%%%%%%%%%%%%%%%%%%%%%%%%%%%%%%%%%%%%%%%%%%%%%%%%%%%%%%%%%%%%%%%%%%%%%%%%%%%%%%%%%%%%%%%%%%%%%%%%%%%%%%%%%%%%%%%%%%%%%%%%%%%%%%%%%%%%%%%%%%%%%%%%%%%%%%%%%%%%%%%%%%%%%%%%%%%%%%%%%%%%%%%%%%%%%%%%%%%%%%%%%%%%%%%%%%%%%%%%%%%%%%%%%%%%%%%%%%%%%

Finally, we will need the following result whose proof can be found in Ladyzhenskaya, \cite[p. $341$]{Ladyzhenskaya}.
\begin{proposition}\label{p2}
%{\bf Mudar o texto. Referencia??}
Assume that $f \in L^p(Q)$ and $u_0 \in W^{2-2/p}_p(\Omega)$. Then the linear problem
\begin{eqnarray*}\label{problemprin00}
\left\{
\begin{array}{ll}
\displaystyle u_t -\Delta u =f(x,t) &\textup{in } Q:=\Omega \times (0,T),\\
\\
\displaystyle  \partial u/\partial \nu = 0 &\textup{on } \sum = \partial \Omega \times (0,T),\\
\\
\displaystyle u(x,0)=u_0(x)  & \textup{for all } x\in \Omega,
\end{array}
\right.
\end{eqnarray*}
has a unique solution $u \in W^{2,1}_p(Q)$ satisfying the estimate
\begin{equation}\label{aux.22}
	\|u\|_{W^{2,1}_p(Q)}\leq C\left(\|u_0\|_{W^{2-2/p}_p(\Omega)} +\|f\|_{L^p(Q)}\right),
\end{equation}
\end{proposition}

%%%%%%%%%%%%%%%%%%%%%%%%%%%%%%%%%%%%%%%%%%%%%%%%%%%%%%%%%%%%%%%%%%%%%%
%%%%%%%%%%%%%%%%%%%%%%%%%%%%%%%%%%%%%%%%%%%%%%%%%%%%%%%%%%%%%%%%%%%%%%
\section{An auxiliary problem}

$\indent$
We  will need  results concerning the following auxiliary nonlinear parabolic boundary value problem:
\begin{equation}\label{P.3}
\left\{
\begin{array}{lll}
\displaystyle \phi_t - \Delta \phi = F(x,t,\phi) + g(x,t) &\textup{if}& (x,t) \in Q,
\\
\displaystyle  \partial \phi/\partial \nu = 0 &\textup{on}& \sum = \partial \Omega \times (0,T),
\\
\displaystyle \phi(x,0)=\phi_0(x) &x \in \Omega. &
\end{array}
\right.
\end{equation}

Let $q$ give by
\begin{equation}\label{q}
q:= \left\{
\begin{array}{lll}
\textup{any number}\ \geq pr &\textup{if}& p\geq \frac{N +2}{2}.
\\
\displaystyle \textup{any number in} \ \left[pr, \frac{p(N+2)}{N+2-2p}\right) &\textup{if}& p <\frac{N +2}{2}.
\end{array}
\right.
\end{equation}
Notice that~(\ref{q}) makes sense by $(H_3)$.

\vspace{0.1cm}
For this auxiliary problem, we have the following existence and regularity result.
\begin{proposition}\label{existpraux1}
Assume that conditions $(H_0)-(H_4)$ hold. Then problem~(\ref{P.3}) admits a unique solution
$\phi \in W^{2,1}_p(Q)$ satisfying the estimate
\begin{equation}\label{estimw21p}
	\|\phi\|_{W^{2,1}_p(Q)}\leq C(1+\|\phi_0\|_{W^{2-2/p}_p(\Omega)} + \|g\|_{L^p(Q)}),
\end{equation}
where $C$ is a constant which depends only $|\Omega|, T, p, r, a, a_0, a_1, a_2, b_0$.

If $\phi_1, \phi_2$ are solutions of~(\ref{P.3}) corresponding to $\phi_{01}, \phi_{02} \in W^{2-2/p}_p(\Omega)$ (instead of $\phi_0$) and $g_1, g_2$ (in place $g$), respectively, with $\|\phi_1\|_{W^{2,1}_p(Q)},\|\phi_2\|_{W^{2,1}_p(Q)}\leq M$, then
\begin{equation}\label{bemposto}
	\|\phi_1 - \phi_2\|_{W^{2,1}_p(Q)}\leq C\left(\|\phi_{01} - \phi_{02}\|_{W^{2-2/p}_p(\Omega)} + \|g_1-g_2\|_{L^p(Q)}\right),
\end{equation}
where $C$ is a constant which depends only $|\Omega|, T, M, p, r, a, a_0, a_1, a_2, b_0$.
\end{proposition}

\subsection{Preparatory results}

To prove Proposition \ref{existpraux1}, we will apply the Leray-Schauder's fixed point theorem.
For this, let us define the nonlinear operator
\begin{equation}
\label{defL1}
\begin{array}{rccl}
L: & L^{pr}(Q)\times [0,1] & \rightarrow & L^{pr}(Q)
\\
&(w,\lambda) & \rightarrow & L(w,\lambda) = \phi,	
\end{array}
\end{equation}
where $\phi$ is the solution of the linear problem
\begin{equation}\label{P.4}
\left\{
\begin{array}{ll}
\displaystyle \phi_t - \Delta \phi = \lambda(F(x,t,w) + g(x,t)) &\textup{in } Q,\\
\\
\displaystyle  \partial\phi/\partial\nu =0 &\textup{on } \sum,\\
\\
\displaystyle \phi(x,0)=\phi_0(x) & \mbox{for all } x \in \Omega.
\end{array}
\right.
\end{equation}

\vspace{0.1cm}
First of all, we have to check that $L$ is well defined. 
In fact, according to (\ref{estF}) we have $F(\cdot,\cdot,w(\cdot)) \in L^{p}(Q)$, $\forall w \in L^{pr}(Q)$.
Then $F(\cdot,\cdot,w(\cdot)) + g \in L^p(Q)$.
Since $\lambda(F(\cdot,\cdot,w(\cdot)) + g) \in L^p(Q)$, by the $L^p$-theory (see Ladyzhenskaya, \cite[pp.~341]{Ladyzhenskaya}),
we have that there exists a unique solution $\phi\in W^{2,1}_p(Q)$ of problem~(\ref{P.4}). It follows from~(\textbf{$H_3$}), (\ref{q}) and Lemma~\ref{imersao continua em Lp 01} that	$W^{2,1}_p(Q) \hookrightarrow L^{pr}(Q)$, and thus $L$ is well defined.

\vspace{0.1cm}
Next, we prove that $L$ has suitable properties.
\begin{lemma}
\label{lemaaux0}
The mapping $L:L^{pr}(Q)\times [0,1] \rightarrow L^{pr}(Q)$ defined in (\ref{defL1}) has the following properties:
\begin{description}
	\item[(i)] $L(\cdot,\lambda):L^{pr}(Q) \rightarrow L^{pr}(Q)$ is compact for every $\lambda \in [0,1]$, i.e., it is continuous and maps bounded sets into relatively compacts sets.
	 \item[(ii)] For every $\epsilon>0$ and every bounded set $A \subset L^{pr}(Q)$ there exists $\delta>0$ such that
	 $$\|L(w,\lambda_1)-L(w,\lambda_2)\|_{L^{pr}(Q)}<\epsilon,$$
	 whenever $w\in A$ and $|\lambda_1 - \lambda_2|<\delta$.
\end{description}
\end{lemma}

\noindent\textbf{Proof:} We start with the proof of {\bf (i)}. 
To prove the continuity of $L(\cdot,\lambda)$, let us consider $w_1, w_2 \in L^{pr}(Q)$ and $\phi_1 = L(w_1,\lambda)$, $\phi_2 = L(w_2,\lambda)$ the corresponding solutions. From~(\ref{defL1}) and~(\ref{P.4}), we obtain
$$
\left\{
\begin{array}{ll}
(\phi_1 - \phi_2)_{t} - \Delta (\phi_1- \phi_2) = \lambda(F(x,t,w_1) - F(x,t,w_2)) & \textup{in } Q,\\
\\
\partial(\phi_1 - \phi_2)/\partial\nu= 0 &\textup{on } \sum,\\
\\
(\phi_1 - \phi_2)(x,0)= 0 & \mbox{for all } x \in \Omega.
\end{array}
\right.
$$

If $w$ is in ${L^{pr}(Q)}$, from Lemma \ref{estF2}, $F(\cdot,\cdot,w)$ is in $L^p(Q)$.
As above, the $L^p$-theory implies the estimate
$$\|\phi_1 - \phi_2\|_{W^{2,1}_p(Q)} \leq C\|F(x,t,w_1) - F(x,t,w_2)\|_{L^p(Q)}.$$

Thus the operator $L(\cdot, \lambda): L^{pr}(Q) \rightarrow W^{2,1}_p(Q)$, given by~(\ref{defL1}), is
continuous. By the compact embedding $W^{2,1}_p(Q) \hookrightarrow L^{pr}(Q)$ (see Lions \cite[pp.~21]{Lions_1983}),
we conclude that $ L(\cdot,\lambda): L^{pr}(Q) 	\rightarrow L^{pr}(Q)$ is continuous and compact for each fixed $\lambda \in [0, 1]$.

\vspace{0.2cm}
Next, we prove {\bf (ii)}.
Let us fix a bounded set $A \subset L^{pr}(Q)$ and consider $(\phi_1, \lambda_1)$, $(\phi_2, \lambda_2)\in L^{pr}(Q)\times[0,1]$
the corresponding solutions of (\ref{P.4}), where we take any $w \in A$. We have From (\ref{defL1}) and (\ref{P.4})
$$
\left\{
\begin{array}{ll}
 (\phi_1 - \phi_2)_t - \Delta (\phi_1- \phi_2) = ( \lambda_1 - \lambda_2)(F(x,t,w) + g(x,t)) &\textup{in } Q,\\
\\
  \partial(\phi_1 - \phi_2)/\partial\nu = 0 &\textup{on } \sum ,\\
\\
 (\phi_1 - \phi_2)(x,0)= 0 & \mbox{for all } x \in \Omega.
\end{array}
\right.
$$

The $L^q$-theory (see Ladyzhenskaya, [\cite{Ladyzhenskaya}, pp.~341]) provides the estimate
$$\|\phi_1 - \phi_2\|_{W^{2,1}_p(Q)} \leq C|\lambda_1 - \lambda_2|(\|F(x,t,w)\|_{L^p(Q)} + \|g\|_{L^p(Q)}).$$

From Lemma \ref{estF2} and (\ref{q}), we have that
$$\|\phi_1 - \phi_2\|_{W^{2,1}_p(Q)} \leq C(A, \|g\|_{L^p(Q)})|\lambda_1 - \lambda_2|,$$
where $C(A, \|g\|_{L^p(Q)})$ is a bounded constant because $A$ is a bounded set.

The proof of \textbf{(ii)} then follows from Lemma \ref{imersao continua em Lp 01}, (\ref{q}) and the last inequality  by taking
$\displaystyle \delta = \epsilon/2C(A, \|g\|_{L^p(Q)})$ for each fixed $\epsilon > 0$.\hfill\caixa

\

In the next result, we prove estimates for any possible fixed points of $L(\cdot,\lambda)$.

\begin{lemma}\label{estapriori0}
Suppose that assumptions $(H_0) - (H_4)$ are satisfied. Then, there exists a number $\rho > 0$, such that
any fixed point $\phi\in L^{pr}(Q)$ of $T(\cdot, \lambda)$ for any $\lambda\in[0,1]$, i. e.,
$T(\phi, \lambda)=\phi$ for some $\lambda\in[0,1]$, satisfies
\begin{equation}\label{aux.13}
	\|\phi\|_{L^{pr}(Q)} < \rho;
\end{equation}
\end{lemma}

\noindent
{\bf Proof:}
Let $\phi\in L^{pr}(Q)$ a fixed point of $T(\cdot, \lambda)$ for some $\lambda\in[0,1]$.
Then, $\phi$ solves the problem
\begin{equation}
\label{P.5}
\left\{
\begin{array}{ll}
\displaystyle  \phi_t - \Delta \phi = \lambda(F(x,t,\phi) + g(x,t)) & \textup{in } Q,\\
\\
\displaystyle  \partial\phi/\partial\nu = 0 & \textup{on } \sum,\\
\\
\displaystyle \phi(x,0)=\phi_0(x) & \mbox{for all } x \in \Omega.
\end{array}
\right.
\end{equation}

Denote
$$Q_{t} =\Omega \times (0,t), \ t \in (0,T].$$

By multiplying first equation in~(\ref{P.3}) by $\phi$, integrating over $Q_t$, using Lemma~\ref{rel.h-F}, Green's formula %, Fubini's theorem
and Young's inequality, we obtain
\[
\begin{array}{c}
\displaystyle
\frac{1}{2}\int_{\Omega}\phi^2 (t)dx + \int_{Q_t}|\nabla \phi|^2dxds
\\
\displaystyle
\leq \frac{1}{2}\left(\|\phi_0\|^2_{L^2(\Omega)} + \|g\|^2_{L^2(Q)}\right) + \frac{1}{2}\int_{Q_t}\phi^2dxds + \int_{Q_t}F(x,s,\phi)\phi dxds.
\end{array}
\]

In view of (\ref{2.2.1}), we get further
$$ \frac{1}{2}\int_{\Omega}\phi^2 (t) dx + \int_{Q_t}|\nabla \phi|^2dxds
\leq C_1\left(1 + \|\phi_0\|^2_{L^p(\Omega)} + \|g\|^2_{L^p(Q)}\right) + C_2\int_{Q_t}\phi^2dxds$$
where $C_1$ and $C_2$ denote positive constants. By Gronwall's inequality, we arrive at

\begin{equation}
\label{desi4}
\displaystyle \frac{1}{2}\int_{\Omega}\phi^2 (t) dx + \int_{Q_t}|\nabla \phi|^2dxds \leq
C_0\left(1 + \|\phi_0\|^2_{W^{2-2/p}_p(\Omega)} + \|g\|^2_{L^p(Q)}\right), \ \forall \ t \in (0,T].	
\end{equation}

As $F(\cdot,\cdot,\phi) + g \in L^p(Q)$. By the $L^p$-theory (see Ladyzhenskaya, \cite[pp.~341]{Ladyzhenskaya},) we have that
\begin{eqnarray*}\label{}
	\|\phi\|_{W^{2,1}_p(Q)}\leq C\left(1+\|\phi_0\|_{W^{2- \frac{2}{p}}_p(\Omega)} +\|F(\cdot,\cdot,\phi)\|_{L^p(Q)}+ \|g\|_{L^p(Q)}\right).
\end{eqnarray*}

By using $q\geq pr$ and the convexity of $|\cdot|^{p/r}$ and Lemma~\ref{estF2}, we get
$$\|F(\cdot,\cdot,\phi)\|_{L^p(Q)} \leq C\|F(\cdot,\cdot,\phi)\|_{L^{q/r}(Q)} \leq C(1+\|\phi\|_{L^{q}(Q)})$$

Next, by (\ref{q}) and Lemma \ref{imersao continua em Lp 01}, we have the sequence of imbeddings
$$\xymatrix{\displaystyle
W^{2,1}_p(Q) \ar @{^{(}->>}[r] & L^q(Q) }  \hookrightarrow L^2(Q),$$
where the first immersion is compact by Lions \cite[pp.~21, Remark~2.3]{Lions_1983}.
The interpolation inequality yields that for all $\epsilon >0$, there exists $C(\epsilon) >0$ such that
\begin{eqnarray*}\label{}
\|v\|_{L^q(Q)} \leq \epsilon\|v\|_{W^{2,1}_p(Q)}+ C(\epsilon)\|v\|_{L^2(Q)}, \ \ \forall v \in W^{2,1}_p(Q).
\end{eqnarray*}

From (\ref{desi4}) and of the two last inequalities, we derive that
 \begin{eqnarray*}\label{}
(1 - \epsilon C)\|\phi\|_{W^{2,1}_p(Q)} \leq
C\left(1+\|\phi_0\|_{W^{2- \frac{2}{p}}_p(\Omega)} + \|\phi_0\|_{W^1_2(\Omega)} + \|g\|_{L^p(Q)}\right),
\end{eqnarray*}
where $\epsilon>0$ is small enough such that $1 - \epsilon C > 0$. So we get
\begin{equation}\label{aux.16}
\begin{array}{lll}
\|\phi\|_{W^{2,1}_p(Q)} &\leq& C\left(1+\|\phi_0\|_{W^{2- \frac{2}{p}}_p(\Omega)} + \|g\|_{L^p(Q)}\right).
\end{array}
\end{equation}

From~(\ref{q}), Lemma~\ref{imersao continua em Lp 01} and~(\ref{aux.16}), we have
\begin{eqnarray}\label{lppr2}
\|\phi\|_{{L^{pr}(Q)}}\leq C\left(1 + \|\phi_0\|_{W^{2- \frac{2}{p}}_p(\Omega)} + \|g\|_{L^{p}(Q)}\right),
\end{eqnarray}

\noindent
which leads to the conclusion that the claim in~(\ref{aux.13}) holds true. \hfill\caixa

%%%%%%%%%%%%%%%%%%%%%%%%%%%%%%%%%%%%%%%%%%%%%%%%%%%%%%%%%%%%%%%%%%%%%%%%%%%%%%%%%%%%%%%%%%%%%%%%%%%%%%%%%%%%%%%%%%%%%%%%%%%%%%%%%%%%%%%%%%%%%%%%%%%%%%%%%%%%%%%%%%%%%%%%%%%%%%%%%%%%%%%%%%%%%%%%%%%%%%%%%%%%%%%%%%%%%%%%%%%%%%%%%%%%%%%%%%%%%%%%%%%%%%%%%%

\subsection{Proof of Proposition \ref{existpraux1}}
From Lemma~(\ref{estapriori0}), we know the existence of a number $\rho > 0$ which satisfies the property stated in~(\ref{aux.13}) and we have that the linear heat equation, $\phi - T(\phi, 0) = 0$, admits a unique solution, i.e., the unique solution of
\[
\left\{ \begin{array}{ll} \displaystyle \phi_t - \Delta \phi =
0 & \;
\; \; \; \; \; \;
\mbox{in }  Q, \\
\displaystyle \partial\phi/\partial\nu = 0 & \; \; \; \; \; \; \;
\mbox{on }  \sum, \\
\phi(x, 0) = \phi_0(x) & \; \; \; \; \; \; \; \mbox{for all } x\in\Omega.
\end{array}\right.
\]

It follows from Leray-Schauder's fixed point theorem (see Friedman \cite[pp.~189, Theorem 3]{Friedman}) that the problem
\begin{equation}\label{L1}
\phi - T(\phi, 1) = 0,
\end{equation}
has a solution $\phi$.

It follows from~(\ref{defL1}) that $\phi$ is a solution of the problem~(\ref{P.3}); from Lemma~\ref{estapriori0} and (\ref{aux.16}) we have $\phi\in W^{2,1}_p(Q)$ and
$$
\|\phi\|_{W^{2,1}_p(Q)}\leq C_1\left(1+ \|\phi_0\|_{W^{2- \frac{2}{p}}_p(\Omega)} + \|g\|_{L^p(Q)}\right).
$$

Since by hypothesis, $\phi_1, \phi_2 \in W^{2,1}_p(Q)$ solve (\ref{P.3}), respectively with $\phi_{01}$, $\phi_{02}$ and $g_1$, $g_2$ , $\phi_1 - \phi_2 \in W^{2,1}_p(Q)$ and  satisfies
\begin{equation}\label{eq.dife}
\left\{
\begin{array}{ll}
(\phi_1 - \phi_2)_t - \Delta (\phi_1- \phi_2) = F(x,t,\phi_1) - F(x,t,\phi_2) + g_1 - g_2 &\textup{in } Q,\\
\\
\partial(\phi_1-\phi_2)/\partial\nu = 0 &\textup{on } \sum,\\
\\
(\phi_1 - \phi_2)(x,0)= \phi_{01} - \phi_{02} & \mbox{for all } x \in \Omega.
\end{array}
\right.
\end{equation}

By multiplying the first equation in~(\ref{P.3}) by $|\phi_1 - \phi_2|^{p-2}(\phi_1 - \phi_2)$, integrating over $Q_t$, and using Green's theorem, we obtain

\[
\begin{array}{c}\label{}
\displaystyle
\frac{1}{p}\int_{\Omega}|\phi_1 - \phi_2|^{p}(t)dx + (p-1)\int_{Q_t}|\nabla (\phi_1 - \phi_2)|^{2}|\phi_1 - \phi_2|^{p-2}dxds
\\
\displaystyle
\leq \frac{1}{p}\int_{\Omega}|\phi_{01} - \phi_{02}|^{p}dx + \frac{2^p}{p}\int_{Q}|g_{1} - g_{2}|^{p}dxds + \frac{p-1}{p}\frac{1}{2^{\frac{p}{p-1}}}\int_{Q_t}|\phi_{1} - \phi_{2}|^{p}dxds
\\
\displaystyle
+ \int_{Q_t}\left(F(x,s,\phi_1) - F(x,s,\phi_2)\right)|\phi_1 - \phi_2|^{p-2}(\phi_1 - \phi_2)dxds.
\end{array}
\]

Due to assumption ($H_1$) and  Gronwall's inequality, it results

\begin{eqnarray}\label{gr1}
\|\phi_1 - \phi_2\|_{L^p(\Omega)}^p\leq C(T,r,p)\left(\|\phi_{01} - \phi_{02}\|_{L^p(\Omega)}^p + \|g_1 - g_2\|_{L^p(Q)}^p\right)
\end{eqnarray}

According to (\ref{q}) and Lemma \ref{imersao continua em Lp 01}, we have $\phi_1, \phi_2 \in L^{pr}(Q)$, which in conjunction with relation (\ref{estF}) yields that $F(\cdot,\cdot,\phi_1) - F(\cdot,\cdot,\phi_2) \in L^p(Q)$. By applying $L_p-$theory to the problem (\ref{eq.dife}) and by using~$(H_0)$, we get the estimate
\begin{equation}
\label{gr2}
\begin{array}{l}
\displaystyle
\|\phi_1 - \phi_2\|_{W^{2,1}_p(Q)}^p\leq C(|\Omega|,T,N)
\\
\displaystyle
\times \left(\|\phi_{01} - \phi_{02}\|_{W^{2-2/p}_p(\Omega)}^p +\|F(\cdot,\cdot,\phi_1) - F(\cdot,\cdot,\phi_2)\|_{L^p(Q)}^p + \|g_1 - g_2\|_{L^p(Q)}^p\right)
\end{array}
\end{equation}

  The inequality $q > pr$ allow us to fix a number $m$ such that
\begin{eqnarray}\label{gr3}
2\leq p\leq \frac{q p}{q + p - pr}\leq pr < m\leq q.
\end{eqnarray}

  Consequently, the next sequence of embeddings holds
\begin{eqnarray}\label{gr4}
W^{2,1}_p(Q) \subset L^{q}(Q) \subset L^{m}(Q) \subset L^{pr}(Q) \subset L^{p}(Q) \subset L^{2}(Q).
\end{eqnarray}

From $(H_2)$, (\ref{gr3}) and H$\ddot{o}$lder's inequality, with $\frac{1}{q_1} + \frac{1}{q_2} =1$, where $q_1 = p/m, q_2 = m/(m-p)$, we conclude that
\begin{equation}
\label{gr5}
\begin{array}{l}
\displaystyle
\|F(\cdot,\cdot,\phi_1) - F(\cdot,\cdot,\phi_2)\|_{L^p(Q)} \leq \left\| G(x,t,\phi_1,\phi_2)^{1/2}|\phi_1 - \phi_2| \right\|_{L^p(Q)}
\\
\hspace{3cm}
\displaystyle
 = \left(\int_QG(x,t,\phi_1,\phi_2)^{p/2}|\phi_1 - \phi_2|^pdxdt\right)^{\frac{1}{p}}
 \\
 \hspace{3cm}
 \displaystyle
\leq \left(\int_QG(x,t,\phi_1,\phi_2)^{n_0/2}dxdt\right)^{\frac{1}{n_0}}\|\phi_1 - \phi_2\|_{L^m(Q)},
\end{array}
\end{equation}

\noindent
where $n_0=mp/(m-p)$. We remark that the previous computations make sense because $G(x,t,\phi_1,\phi_2)^{n_0/2} \in L^1(Q)$. Indeed, taking into account the growth condition in $(H_2)$, $G(x,t,\phi_1,\phi_2) \in L^{\frac{q}{2(r-1)}}(Q)$, whenever $\phi_1, \phi_2 \in L^{q}(Q)$, and by (\ref{gr3}) it is true that
\begin{eqnarray}\label{gr6}
\frac{q}{r-1} > n_0 > 2.
\end{eqnarray}

Then, the inequalities in~(\ref{gr6}) lead to the claim above. Combining~(\ref{gr2}) and~(\ref{gr5}), we arrive at
\begin{equation}
\label{gr7}
\begin{array}{l}
\displaystyle
\|\phi_1 - \phi_2\|_{W^{2,1}_p(Q)}^p\leq C(|\Omega|,T,N)\left( \|\phi_{01} - \phi_{02}\|_{W^{2-2/p}_p(\Omega)}^p\right.
\\
\hspace{2cm}
\displaystyle
\left. +\|G(x,t,\phi_1,\phi_2)\|_{L^{\frac{n_0}{2}}(Q)}^{1/2}\|\phi_1 - \phi_2\|_{L^m(Q)} + \|g_1 - g_2\|_{L^p(Q)}^p\right)
\end{array}
\end{equation}
In addition, we have for any  $(x,t) \in Q$ that
\[
\begin{array}{l}
\displaystyle
\left(1 + |\phi_1(x,t)|^{2r-2} + |\phi_2(x,t)|^{2r-2}\right)^{n_0/2}
\\
\hspace{2cm}
\displaystyle
\leq C(p,r) (1 + |\phi_1(x,t)|^{n_0(r-1)} + |\phi_2(x,t)|^{n_0(r-1)} ).
\end{array}
\]

This last relation, inequality~(\ref{gr6}) and estimate~(\ref{gr7}) then imply
\begin{equation}
\label{gr8}
\begin{array}{l}
\displaystyle
\|\phi_1 - \phi_2\|_{W^{2,1}_p(Q)}
\\
\displaystyle
\leq C(|\Omega|,T,N) \left[ \|\phi_{01} - \phi_{02}\|_{W^{2-2/p}_p(\Omega)} + \|g_1 - g_2\|_{L^p(Q)} \right.
\\
\hspace{0.5cm}
\displaystyle
\left. +C(|\Omega|, T, r, p, b_0)\left( 1 + \|\phi_1\|_{L^{n_0(r-1)}}^{r-1} + \|\phi_2\|_{L^{n_0(r-1)}}^{r-1}\right)\|\phi_1 - \phi_2\|_{L^m(Q)} \right]
\vspace{0.1cm}
\\
\displaystyle
\leq C(|\Omega|,T,N,p,r,b_0) \left[ \|\phi_{01} - \phi_{02}\|_{W^{2-2/p}_p(\Omega)} + \|g_1 - g_2\|_{L^p(Q)} \right.
\\
\hspace{0.5cm}
\displaystyle
\left. +\left( 1 + \|\phi_1\|_{L^{n_0(r-1)}}^{r-1} + \|\phi_2\|_{L^{n_0(r-1)}}^{r-1}\right)\|\phi_1 - \phi_2\|_{L^m(Q)} \right]
\vspace{0.1cm}
\\
\displaystyle
\leq C(|\Omega|,T,N,p,r,b_0)(1 + 2M^{r-1})
\\
\hspace{0.5cm}
\displaystyle
\times \left[ \|\phi_{01} - \phi_{02}\|_{W^{2-2/p}_p(\Omega)} + \|g_1 - g_2\|_{L^p(Q)} +\|\phi_1 - \phi_2\|_{L^m(Q)} \right]
\end{array}
\end{equation}

Due to the embeddings in~(\ref{gr4}), the interpolation inequality (see Lions~\cite[pp.~58, Lemma~5.1]{Lions}) yields that for all
$\epsilon>0$, there exists $C(\epsilon)>0$ such that
\begin{eqnarray}\label{gr9}
\|\phi\|_{L^m(Q)}  \leq \epsilon \|\phi\|_{W^{2,1}_p(Q)} +  \|\phi\|_{L^p(Q)}, \forall \phi \in W^{2,1}_p(Q).
\end{eqnarray}

Thus, from (\ref{gr1}), (\ref{gr8}) and (\ref{gr9}), we get that
\begin{equation}
\label{gr10}
\begin{array}{c}
\displaystyle
(1 - \epsilon C(|\Omega|,T,N,p,r,b_0))\|\phi_1 - \phi_2\|_{W^{2,1}_p(Q)}
\\
\displaystyle
\leq C(|\Omega|,T,N,p,r,a_0,b_0) \left[ \|\phi_{01} - \phi_{02}\|_{W^{2-2/p}_p(\Omega)} + \|g_1 - g_2\|_{L^p(Q)} \right].
\end{array}
\end{equation}

Next, by taking $\epsilon >0$ such that $1 - \epsilon C(|\Omega|,T,N,p,r,b_0) >0$, the last inequality implies the estimate~(\ref{bemposto}),
completing the proof of Proposition~\ref{existpraux1}. \hfill\caixa

%%%%%%%%%%%%%%%%%%%%%%%%%%%%%%%%%%%%%%%%%%%%%%%%%%%%%%%%%%%%%%%%%%%%%%%%%%%%%%%%%%%%%%%%%%%%%%%%%%%%%%%%%%%%%%%%%%%%%%%%%%%%%
\section{Proof of the main result}

\subsection{Preparatory results}
To prove Theorem \ref{meantheorem}, we will again apply the Leray-Schauder's fixed point theorem. 
Consider the nonlinear operator $\mathcal{L}$ defined by
\begin{equation}
\label{defL}
\begin{array}{rccl}
\mathcal{L}: & L^p(Q)\times [0, 1] & \rightarrow & L^p(Q)
\\
& (g,\lambda) & \rightarrow & \mathcal{L}(g,\lambda) = u,
\end{array}
\end{equation}

\noindent
where $u$ is the unique solution of the following linear parabolic boundary value problem
\begin{equation}\label{P.6}
\left\{ \begin{array}{ll}
\displaystyle
u_t -\Delta u =\lambda(-l\phi_t + f) & \mbox{ in } Q,
\\
\displaystyle
\partial u/\partial\nu =0 & \mbox{ on} \sum,
\\
\displaystyle
u(x, 0) = \lambda{u}_0(x) & \mbox{ for all } x\in\Omega,
\end{array}\right.
\end{equation}
with $\phi$ the unique solution of the following semilinear parabolic boundary value problem (see Proposition~\ref{existpraux1})
\begin{equation}
\label{P.6.1}
\left\{ \begin{array}{ll}
\displaystyle
\phi_t -\Delta\phi =F(x, t, \phi) +g & \mbox{ in } Q,
\\
\displaystyle
\partial\phi/\partial\nu =0 & \mbox{ on }  \sum,
\\
\displaystyle
\phi(x, 0) = {\varphi}_0(x) & \mbox{ for all } x\in\Omega.
\end{array}\right.
\end{equation}

First of all, we observe that the operator $\mathcal{L}$ is well defined. 
In fact, according to Proposition \ref{existpraux1}, since $g\in L^p(Q)$, there is only one solution $\phi \in W^{2,1}_p(Q)$ of~(\ref{P.6.1}). And since, $-l\phi_t + f \in L^p(Q)$, it follows from the Proposition~\ref{p2}, that~(\ref{P.6}) has a unique solution $u \in W^{2,1}_p(Q)$.
Therefore, $\mathcal{L}(g,\lambda) = u \in W^{2,1}_p(Q) \subset L^p(Q)$, 
and $\mathcal{L}$ is well defined.

\begin{lemma}
\label{lemaaux1}
The mapping $\mathcal{L}:L^p(Q)\times [0,1] \rightarrow L^p(Q)$ has the following properties:
\begin{description}
	\item[(i)] $\mathcal{L}(\cdot,\lambda):L^p(Q) \rightarrow L^p(Q)$ is compact for every $\lambda \in [0,1]$, i. e.,
it is con\-tin\-uous and maps bounded sets into relatively compacts sets.
	 \item[(ii)] For every $\epsilon>0$ and every bounded set $A \subset L^p(Q)$ there exists $\delta>0$ such that,
	 whenever $g \in A$ and $|\lambda_1 - \lambda_2|<\delta$, there holds
	 \[
	 \|\mathcal{L}(g,\lambda_1)- \mathcal{L}(g,\lambda_2)\|_{L^p(Q)}<\epsilon.
	 \]
\end{description}
\end{lemma}

\noindent\textbf{Proof:} We start by proving {\bf (i)}.
To prove the continuity of $\mathcal{L}(\cdot,\lambda)$, let us consider $g_1$, $g_2 \in L^p(Q)$ and $u_1 = \mathcal{L}(g_1,\lambda)$,
$u_2 = \mathcal{L}(g_2,\lambda)$ the corresponding solutions. By (\ref{defL}) and (\ref{P.6}) we derive

\[  %begin{equation}\label{P.7}
\left\{
\begin{array}{ll}
\displaystyle
(u_1 - u_2)_t - \Delta(u_1 - u_2) = -\lambda l(\phi_1 - \phi_2)_t &  \mbox{ in }  Q,
\\
\displaystyle
\partial(u_1 - u_2)/\partial\nu = 0 & \mbox{ on } \sum,
\\
\displaystyle
(u_1 - u_2)(x, 0) = 0 & \mbox{ for all } x\in\Omega.
\end{array}\right.
\]  %end{eqnarray}

We can then apply the classical $L^q$-theory (see Ladyzhenskaya, \cite[pp.~341]{Ladyzhenskaya}),
and obtain the estimate
\begin{eqnarray}\label{n3.3}
\|u_1 - u_2\|_{W^{2,1}_p(Q)} \leq C\|(\phi_1 - \phi_2)_t\|_{L^p(Q)}.
\end{eqnarray}

Estimate (\ref{bemposto}) in Proposition \ref{existpraux1} applied to (\ref{P.6.1}) choosing $\phi_{01}=\phi_{02}=\phi_{0}$ insures that
\begin{eqnarray}\label{n3.4}
\|\phi_1 - \phi_2\|_{W^{2,1}_p(Q)} \leq C\|g_1 -g_2\|_{L^p(Q)}.
\end{eqnarray}

From (\ref{n3.3}) and (\ref{n3.4}) yield a new constant $C > 0$ such that
\begin{eqnarray}\label{n3.5}
\|u_1 - u_2\|_{L^p(Q)} \leq C\|g_1 -g_2\|_{L^p(Q)},
\end{eqnarray}

\noindent
which implies the continuity of $\mathcal{L}(\cdot,\lambda)$.

Furthermore, $\mathcal{L}(\cdot,\lambda)$ is a compact operator. 
In fact,  from (\ref{n3.5}) we have that $\mathcal{L}(\cdot,\lambda)$, 
regarded as a mapping from the space $L^p(Q)$ into the space $W^{2,1}_p(Q)$,  is continuous. 
From the compactness of the embedding $W^{2,1}_p(Q) \hookrightarrow L^p(Q)$ (see Lions \cite[pp.~21, Remark~2.3]{Lions_1983}), 
we conclude that $\mathcal{L}(\cdot,\lambda): L^p(Q)\rightarrow L^p(Q)$ is continuous and compact for each fixed $\lambda \in [0, 1]$.

\vspace{0.2cm}
Next we prove {\bf (ii)}. For this, fix a bounded set $A \subset L^p(Q)$;
consider $\lambda_1$, $\lambda_2\in [0,1]$ and $g\in A$,
and set $u_1~=~\mathcal{L}(g, \lambda_1)$, $u_2 = \mathcal{L}(g, \lambda_2)$. 
Then we have
\begin{equation}
\label{P.8}
\left\{
\begin{array}{ll}
\displaystyle
(u_1 - u_2)_t -\Delta(u_1 - u_2) =-l(\lambda_1{\phi_1}_t -\lambda_2{\phi_2}_t)
+(\lambda_1 -\lambda_2)f & \mbox{ in } Q,
\\
\displaystyle
\partial(u_1 - u_2)/\partial\nu =0 & \mbox{ on } \sum,
\\
\displaystyle
(u_1 - u_2)(x, 0) =(\lambda_1 - \lambda_2)u_0 &  \mbox{ for all } x\in\Omega.
\end{array}
\right.
\end{equation}

By an argument similar to that used to derive (\ref{n3.4}), we get
\begin{eqnarray}\label{n3.7}
\|\phi_1 - \phi_2\|_{W^{2,1}_p(Q)} \leq C|\lambda_1 -
\lambda_2|\|g\|_{L^p(Q)}\leq C(A)|\lambda_1 -
\lambda_2|.
\end{eqnarray}

By applying the classical $L^q$-theory (see Ladyzhenskaya, \cite[pp.~341]{Ladyzhenskaya}), for the first equation of~(\ref{P.8}) leads to the inequality
\[
\begin{array}{l}
\|u_1 -u_2\|_{W^{2,1}_p(Q)}
\\
\hspace{0.5cm}
\displaystyle
\leq
C\left( \|\lambda_1{\phi_1}_t -\lambda_2{\phi_2}_t\|_{L^p(Q)}
+|\lambda_1 -\lambda_2|\|f\|_{L^p(Q)} +|\lambda_1 -\lambda_2|\|u_0\|_{W^{2-2/p}_p(\Omega)}\right)
\\
\hspace{0.5cm}
\displaystyle
\leq C|\lambda_1 -\lambda_2|\left(\|{\phi_1}_t\|_{L^p(Q)}+ \|u_0\|_{W^{2-2/p}_p(\Omega)} + \|f\|_{L^p(Q)}\right)
+C\lambda_2\|{\phi_1}_t - {\phi_2}_t\|_{L^p(Q)}.
\end{array}
\]

Now, estimate (\ref{estimw21p}) ensures that $\|{\phi_1}_t\|_{L^p(Q)}$ is bounded because $A$ is bounded;
 therefore, from (\ref{n3.7}), we get
\begin{eqnarray*}
\label{n3.8}
\|u_1 - u_2\|_{L^p(Q)} \leq C\|u_1 - u_2\|_{W^{2,1}_p(Q)} \leq C(A)|\lambda_1 -
\lambda_2|.
\end{eqnarray*}
It follows from (\ref{n3.7}) and the last inequality  the proof of \textbf{(ii)} with $\delta =\epsilon/(2 C(A))$, for each fixed $\epsilon > 0$.\hfill\caixa

\vspace{0.3cm}
In the next result, we prove estimates for any possible fixed points of $\mathcal{L}(\cdot,\lambda)$.
\begin{lemma}
\label{estapriori}
Under assumptions $(H_0) - (H_4)$, there exists a number $\rho > 0$, independent of $\lambda \in [0,1]$, with the following property
\begin{equation}\label{aux.24}
	\mathcal{L}(u,\lambda)=u	\textup{, for any } \lambda\in [0,1], \ \Rightarrow \ \|u\|_{L^p(Q)} < \rho.
\end{equation}
\end{lemma}

\noindent
{\bf Proof:}
Let us consider $u$ a possible fixed point of $\mathcal{L}(\cdot,\lambda)$. 
Then $u$ satisfies $(\ref{P.6})$, where $\phi$ is the unique solution of the nonlinear parabolic boundary value problem (see Proposition~\ref{existpraux1})
\[
\left\{
\begin{array}{ll}
\displaystyle
\phi_t -\Delta\phi =F(x,t,\phi) + u & \mbox{ in } Q,
\\
\displaystyle
\partial\phi/\partial\nu =0 & \mbox{ on } \sum,
\\
\displaystyle
\phi(x, 0) = {\phi}_0(x) & \mbox{ for all } x\in\Omega.
\end{array}
\right.
\]

By Proposition~\ref{p2}, $u\in L^p(Q)$. Then, by Proposition~\ref{existpraux1}, there is a unique $\phi\in W^{2,1}_p(Q)$ solving the above problem. In addition, estimate~(\ref{estimw21p}) is true with $u$ replacing $g$. From~(\ref{estimw21p}), we get
\begin{eqnarray}\label{P.10.1}
\| \phi_t\|_{L^p(Q)}\leq \|\phi\|_{W^{2,1}_p(Q)}\leq C\left(1+\|\phi_0\|_{W^{2-2/p}_p(\Omega)} + \|u\|_{L^p(Q)}\right),
\end{eqnarray}
where $C$ is a constant which depends only $|\Omega|, T, p, r, a, c_0, d_0$.

By multiplying the first equation in (\ref{P.6}) by $|u|^{p-2}u$, integrating over $Q_t$, with $t\in(0,T]$, and using %Fubini's theorem,
Green's formula and Young's inequality, we obtain
\[
\begin{array}{c}
\displaystyle
\frac{1}{p}\int_{\Omega}|u|^{p}(t) dx + (p-1)\int_{Q_t}|\nabla u|^{2}|u|^{pr-2}dxds
\\
\displaystyle
\leq \frac{1}{p}\|u_0\|^{p}_{L^p(\Omega)} + \frac{1}{p}\int_{Q_t}|\phi_t|^p dxds + \frac{2(p-1)}{p}\int_{Q_t}|u|^{p}dxds + \frac{1}{p}\int_{Q_t}|f|^{p}dxds,
\end{array}
\]
for all $t \in (0,T]$. By combining (\ref{P.10.1}) (with $Q_t$ in place of $Q$) and the last inequality%(\ref{P.11})
, it turns out that
\[
\begin{array}{cc}
\displaystyle
\frac{1}{p}\int_{\Omega}|u|^{p}(t) dx + (p-1)\int_{Q_t}|\nabla u|^{2}|u|^{pr-2}dxds
\\
\displaystyle
\leq \nonumber \frac{1}{p}\|u_0\|^{p}_{L^p(\Omega)} + \frac{C}{p}\left(1+\|\phi_0\|^{p}_{W^{2-2/p}_p(\Omega)}\right)
\\
\displaystyle+ \left[\frac{C}{p} + \frac{2(p-1)}{p}\right]\int_{Q_t}|u|^{p}dxds + \frac{1}{p}\int_{Q_t}|f|^{p}dxds, & \forall \ t \in (0,T].
\end{array}
\]

By using Gronwall's lemma in last inequality, we obtain a positive constant $C=C(|\Omega|, T, p, r, a, c_0, d_0)$ such that
\begin{eqnarray}
\label{P.13}
\|u\|_{L^p(Q)} \leq C (1 + \|u_0\|_{L^p(\Omega)} + \|\phi_0\|_{W^{2-2/p}_p(\Omega)} + \|f\|_{L^p(Q)} ).
\end{eqnarray}

By taking for instance $\rho = C (1 + \|u_0\|_{L^p(\Omega)} + \|\phi_0\|_{W^{2-2/p}_p(\Omega)} + \|f\|_{L^p(Q)} ) +1$,
we then have property (\ref{aux.24}). \hfill\caixa

%%%%%%%%%%%%%%%%%%%%%%%%%%%%%%%%%%%%%%%%%%%%%%%%%%%%%%%%%%%%%%%%%%%%%%%%%%%%%%%%%%%%%%%%%%%%%%%%%%%%%%%%%%%%%%%%%%%%%%%%%%%%%%%%%%%%%%%%%%%%%%%%%%%%%%%%%%%%%%%%%%%%%%%%%%%%%%%%%%%%%%%%%%%%%%%%%%%%%%%%%%%%%%%%%%%%%%%%%%%%%%%%%%%%%%%%%%%%%%%%%%%%%%%%%%

\subsection{Proof of Theorem \ref{meantheorem}}
From Lemma (\ref{estapriori}), we know the existence of a number $\rho > 0$ which satisfies the property stated in (\ref{aux.24});
moreover, any fixed point $\bar{u}$ corresponding to $\lambda =0$, that is, $\bar{u}$ satisfying $\mathcal{L}(\bar{u},0) = \bar{u}$,  
is exactly the solution of the standard linear heat equation
\[
\left\{ \begin{array}{ll}
\displaystyle
\bar{u}_t -\Delta \bar{u} =0 & \mbox{ in } Q,
\\
\displaystyle
\partial \bar{u}/\partial \nu =0 & \mbox{ on } \sum,
\\
\displaystyle
\bar{u}(x, 0) = 0 & \mbox{ for all } x\in\Omega,
\end{array}\right.
\]

\noindent
which has a unique solution.

It then follows from Leray-Schauder's fixed point theorem (see Friedman \cite[pp.~189, Theorem~3]{Friedman}) that problem
\[
\mathcal{L}(u, 1) = u,
\]
has a solution $u$.

This implies that equation (\ref{P.6}) has a solution $u \in W^{2,1}_p(Q)$, with a corresponding
$\phi \in W^{2,1}_p(Q)$ which is the unique solution of (\ref{P.6.1}) with $g = u$.
That is, there is a solution $(u,\phi) \in W^{2,1}_p(Q)\times W^{2,1}_p(Q)$ of (\ref{P.1}), which, by estimate (\ref{estimw21p}) in Proposition \ref{existpraux1}, satisfies
\begin{equation}\label{P.14}
	\|\phi\|_{W^{2,1}_p(Q)}\leq C(|\Omega|, T, p, r, a, c_0, d_0)
\left(1+\|\phi_0\|_{W^{2-2/p}_p(\Omega)} + \|u\|_{L^p(Q)}\right).
\end{equation}

By using Proposition~\ref{p2} with first equation in (\ref{P.1}), combined with the last inequality, we get
\[
	\begin{array}{l}
\displaystyle
\|u\|_{W^{2,1}_p(Q)}
\leq  C\left(\|u_0\|_{W^{2-2/p}_p(\Omega)} + \| \phi_t\|_{L^p(Q)} + \|f\|_{L^p(Q)}\right)
\\
\displaystyle
\phantom{\|u\|_{W^{2,1}_p(Q)}}
\leq C\left(1 + \|\phi_0\|_{W^{2-2/p}_p(\Omega)} + \|u_0\|_{W^{2-2/p}_p(\Omega)}+ \|u\|_{L^p(Q)} + \|f\|_{L^p(Q)}\right).	
	\end{array}
\]

From this result and (\ref{P.13}), we obtain that there is a positive constant $C=C(|\Omega|,T, p, r, a, c_0, d_0)$
such that
\begin{equation}
\label{first est}
\|u\|_{W^{2,1}_p(Q)} \leq C ( 1 + \|\phi_0\|_{W^{2-2/p}_p(\Omega)} + \|u_0\|_{W^{2-2/p}_p(\Omega)} + \|f\|_{L^p(Q)}).
\end{equation}

By collecting the results in estimates (\ref{P.13}), (\ref{P.14}) and (\ref{first est}), we conclude that there is a positive constant $C=C(|\Omega|,T, p, r, a, c_0, d_0)$ such that
arrive the inequality
\[
\|u\|_{W^{2,1}_p(Q)} + \|\phi\|_{W^{2,1}_p(Q)} \leq  C (1 + \|\phi_0\|_{W^{2-2/p}_p(\Omega)} + \|u_0\|_{W^{2-2/p}_p(\Omega)} + \|f\|_{L^p(Q)}),	
\]
which proves (\ref{estprinc}).

\vspace{0.2cm}
Next, we prove the uniqueness of the solution.

For this, take $(\phi_1, u_1)$ and $(\phi_2, u_2)$ two solutions as in the statement of Theorem \ref{meantheorem} corresponding to the same $\phi_0$, $u_0$ and $f$.
By denoting $u:=u_1 - u_2$ and $\phi:=\phi_1 - \phi_2$, we have that $u$ and $\phi$ satisfy the following problem:
\begin{eqnarray}\label{P.10}
\left\{
\begin{array}{ll}
\displaystyle
u_t -\Delta u =-l\phi_t & \mbox{ in } Q,
\\
\displaystyle
\phi_t -\Delta\phi =F(x,t,\phi_1) -F(x,t,\phi_2) +u & \mbox{ in } Q,
\\
\displaystyle
\partial u/\partial \nu =\partial\phi/\partial\nu =0 & \mbox{ on } \sum,
\\
\displaystyle
u(x,0) =\phi(x,0) =0 & \mbox{ for all } x\in\Omega.
\end{array}
\right.
\end{eqnarray}

To complete the proof, it is enough to show that $\|u\|_{L^2(Q)}=0$ and $\|\phi\|_{L^2(Q)}=0$. 
For this, we multiply first equation in (\ref{P.10}) by $u + l\phi$, integrate on $Q_t$ and use Green's formula and Young's inequality to get
\begin{eqnarray*}
\frac{1}{2}\int_\Omega(u + l\phi)^2 (t) dx +
\int_{Q_t}(|\nabla u|^2 + l\nabla u.\nabla\phi)dxds=0.
\end{eqnarray*}

By multiplying the second equation in  (\ref{P.10}) by $\phi$, integrating on $Q_t$ and using the Green's formula \textbf{$(H_1)$} and Young's inequality, we conclude that
\[
\begin{array}{c}
\displaystyle
\frac{1}{2}\int_\Omega\phi^2 (t) dx +\int_{Q_t}|\nabla\phi|^2dxds
\\
\displaystyle
=\int_{Q_t}(F(x, s, \phi_1) -F(x, s, \phi_2))\phi\ dxds + \int_{Q_t}u.\phi\ dxds
\\
\displaystyle
\leq C\left(\int_{Q_t}\phi^2dxds + \int_{Q_t}u^2dxds\right).
\end{array}
\]

Next, by multiplying this last inequality by an arbitrary positive constant $A$ and adding the result to the previous inequality, we obtain
\[
\begin{array}{c}
\displaystyle
\frac{1}{2}\int_\Omega((u + l\phi)^2 + A\phi^2)(t) dx
+\int_{Q_t}(|\nabla u|^2 + l\nabla u.\nabla\phi +A|\nabla\phi|^2)dxds
\\
\displaystyle
\leq C\int_{Q_t}(u^2 +\phi^2)dxds,
\end{array}
\]

By using suitably H\"older's inequality and taking $A = 1 + l^2$ and then Gronwall's lemma, we arrive at
\begin{eqnarray*}
\int_\Omega(u^2 + \phi^2)(t) dx + \int_{Q_t}(|\nabla u|^2
+ |\nabla\phi|^2)dxds \leq 0.
\end{eqnarray*}

Therefore, $\|u\|_{L^2(Q)}=0$ and $\|\phi\|_{L^2(Q)}=0$ and thus, $\phi_1=\phi_2$ e $u_1=u_2$.
\hfill\caixa

%%%%%%%%%%%%%%%%%%%%%%%%%%%%%%%%%%%%%%%%%%%%%%%%%%%%%%%%%%%%%%%%%%%%%%%%%%%%%%%%%%%%%%%%%%%%%%%%%%%%%%%%%%%%%%%%%%%%%%%%%%%%%%%%%%%%%%%%%%%%%%%%%%%%%%%%%%%%%%%%%%


\begin{thebibliography}{99}




%\bibitem{Adams} Adams, R.A., Sobolev Spaces, Academic Press, New York, 1975.

\bibitem{Bates} P. W. Bates  and S. Zheng, Inertial manifolds and inertial sets for the phase field equations, J Dynam. Differential Equations 4(1992), 375-398.

%\bibitem{Brezis} Brezis, H.. Analyse Fonctionnelle. Th\'eorie et  Applicatons. Masson, Paris. 1983.

\bibitem{Caginalp} G. Caginalp, An analysis of a phase field model of a free boundary. Arch. Rational Mech. Anal. 92 (1986), 205-245.

\bibitem{OAC} O. C\^arja, A. Miranville and C. Moro\c{s}anu, On the existence, uniqueness and regularity of solutions to the phase-field system with a general regular potential and a general class of nonlinear and non-homogeneous boundary conditions, Nonlinear Analysis 113 (2015), 190-208.

%\bibitem{Vespri} Cannarsa, P. and Vespri, V., On maximal $L^p$ regularity for the abstract Cauchy problem. Boll. Un. Mat. Ital. B (6) 5 (1986),  no. 1, 165--175.

%\bibitem{Da Prato and Zabczyk} Da Prato, J. and Zabczyk, J., Stochastic Equations in Infinite Dimensions. Cambridge University Press, 1992.


%\bibitem{Elliot Zheng} Elliot, C. M. and Songmmu Zheng, Global existence and stability of solutions to the phase field equations, in: \textit{Free Boundary Problems}, K.-H. Hoffman and J. Sprekls eds., Birkhauser Verlag, Basel, (1990), 46-58.

%\bibitem{Evans} Evans,  L.C., Partial Differential Equations, American Mathematical Society, 1998.

%\bibitem{Fonseca-Gangbo} Fonseca, I. and Gangbo, W., Degree Theory in Analysis and Applications, Clarendon, Oxford, 1995.

\bibitem{Friedman} A. Friedman, Foundations of Modern Analysis, Hold, Rinehart and Winston, INC, 1970.

%\bibitem{Guidetti} Guidetti, D., On interpolation with boundary conditions, Math. Z., 207(1991), 439-460.

%\bibitem{Henry}  Henry, D., Geometric Theory of Semilinear Parabolic Equations, Lecture notes in Mathematics 840, Springer-Verlag, 1981.


\bibitem{Hoffman} K.H. Hoffman and L. Jiang, Optimal control problem of a phase field model for solidification, Numer. Funct. Anal., 13 (1992), 11-17.


%\bibitem{Yamada} Hoshino, H. and  Yamada, Y., Solvability and smoothing effect for Semilinear Parabolic Equations, Funkt. Ekv., 34 (1991), 475-494.

%\bibitem{Lions-Magenes} Lions, J. L. and Magennes, E., Probl\'emes aux limites non homog\'enes et application II, Dunod, Paris, 1968.

%\bibitem{Lunardi} Lunardi, A., Analytic Semigroups and Optimal Regularity in Parabolic Problems, Birkh$\ddot{a}$user, Basel (1995).

%\bibitem{Pazy}  A. Pazy. Semigroups of Linear Operator and Applications to Partial Differential Equations. Springer-Verlag, 1983.

%\bibitem{Simon} Simon, J.,Compact sets in the space $L^p(0,T;B)$, Annali Mat. Pura Aplli., Serie IV, v.146 (1987), 65-96.

%\bibitem{Taylor} Taylor, M.E., Partial Differential Equations III,
%Nonlinear Equations, Springer Verlag, 1997.

%\bibitem{Kufner} Kufner, A., John, O., Fu$\check{c}$ik, S., Function spaces, Noordhoff Int. Publ., Leyden, 1977.

\bibitem{Ladyzhenskaya} O. Ladyzhenskaya, V. Solonnikov and N. Uraltseva, Linear and Quasilinear Equations
of Parabolic Type, American Mathematical Society, 1968.

\bibitem{Lions} J.L. Lions, Quelques M\'ethodes de R\'esolution des Problemes aux Limites  non Lin\'eaires, Dunod and Guthier-Villars, 1969.

\bibitem{Lions_1983} J.L. Lions, Contr\^ole des Syst\`emes Distribu\'es Singuliers, M\'ethodes Math\'ematiques de I'Informatique, Gautier-Villars, 1983.

\bibitem{Morosanu-Motreanu 2} C. Moro\c{s}anu and D. Motreanu, A Generalized Phase-Field System. Journal of Math. Analysis and Applications 237 (1999), 515-540.

\bibitem{Morosanu-Motreanu 1} C. Moro\c{s}anu and D. Motreanu, The phase field system with a general nonlinearity. Int. J. Differ. Equ. Appl.  1  (2000),  no. 2, 187-204.





%\bibitem{Peetre} Peetre, J., New Thoughts on Besov Spaces, Tekniska Hogskolan Lund, 1976.

%\bibitem{Schwartz} Schwartz, J.T., Nonlinear Functional Analysis. Gordon and Breach, New York, 1969.

%\bibitem{Solonnikov and Garroni} Garroni, M.G. and Solonnikov, V. A., On parabolic oblique derivative problem with Holder continuous coeficients. Comm. in Partial Differential Equations, 9(14), 1323 - 1372(1984).

%\bibitem{Wahl} Von Wahl, W., The equation $u'+A(t)u=f$ in a Hilert space and $L^p$-estimates for parabolic equations. J. London Math. Soc., 25 (1982) 483 - 497.

%\bibitem{Weidemaier} Weidemaier, P., On the Sharp Initial Trace of Functions with Derivates in $L_q(0,T;L_p(\Omega))$. Bolletino U. M. I. (7) 9 - B(1995), 321-338.


\end{thebibliography}
\end{document}